\documentclass{amsart}

\usepackage{amsmath,amssymb}

\newtheorem{theorem}{Theorem}
\newcommand{\bt}{\begin{theorem}}
\newcommand{\et}{\end{theorem}} 

\newtheorem{corollary}{Corollary}
\newcommand{\bc}{\begin{corollary}}
\newcommand{\ec}{\end{corollary}} 

\newtheorem{lemma}{Lemma}
\newcommand{\bl}{\begin{lemma}}
\newcommand{\el}{\end{lemma}}

\newtheorem{problem}{Problem}
\newcommand{\bprob}{\begin{problem}}
\newcommand{\eprob}{\end{problem}}

\newcommand{\beq}{\begin{equation}}
\newcommand{\eeq}{\end{equation}}
\newcommand{\N}{\ensuremath{ \mathbf N }}
\newcommand{\Z}{\ensuremath{\mathbf Z}}

\newcommand{\benum}{\begin{enumerate}}
\newcommand{\eenum}{\end{enumerate}}

\newcommand{\mca}{\mathcal{A}}

\DeclareMathOperator{\qqand}{\qquad\text{and}\qquad}

\DeclareMathOperator{\diam}{\text{diam}}

\title[Bases of order $h$ for $n$]{Problems in additive number theory, VII:  \\
The structure of additive $h$-bases for $n$}

\author{Melvyn B. Nathanson}

\address{Department of Mathematics, Lehman College (CUNY), Bronx, NY, USA}  
\email{melvyn.nathanson@lehman.cuny.edu}

\date{\today}


\begin{document}

\begin{abstract}
In additive number theory, a finite set $A$ of integers is an 
\emph{$h$-basis for $n$} if every integer
in $\{0,1,2,\ldots, n\}$ can be represented as the sum of exactly 
 $h$  not necessarily distinct elements of $A$.
This paper introduces a new class of problems  for these and related additive bases.  
The problems are designed, in part, to be susceptible to solution by AI. 
\end{abstract}

\maketitle

\section{What can humans do?}
In \emph{Magnifica Humanitas}~\cite[Chapter 3]{leo26},  Pope Leo XIV wrote 
\begin{quotation}
 It is not possible to provide a single, comprehensive definition of AI.  
 What can be stated, however, is that we must avoid the misconception of equating 
 this type of “intelligence” with that of human beings. 
 \end{quotation}
Artificial intelligence systems (AI) 
today\footnote{This paper was uploaded to arXiv on May 25, 2026, and published on arXiv on May 26, 2027.  Pope Leo's  encyclical on AI was also released on May 25, 2026.} 
can solve mathematical problems sufficiently difficult to constitute an acceptable  
Ph.D. thesis and they have even solved a few significantly more important problems.  
But AI has not yet been able to decide what is an ``interesting'' problem nor
to formulate new conjectures  or new areas of research.  
That is the (likely transitory) advantage that we humans retain. 

Perhaps the future job description of a mathematician will be:  
\benum
\item
Decide what problems are easy enough for AI to solve.
\item
Invent problems too hard for AI to solve. 
\item
Solve problems too hard for AI to solve.
\eenum

This paper introduces some new questions in additive number theory.  All should be solvable 
by humans and some by the AI systems now being applied to  questions in arithmetic. 
An asterisk (*) identifies problems that AI might have difficulty solving.

\section{Bases of order $h$ for $n$} 
Let $\N= \{1,2,3,\ldots \}$  be the set of positive integers, $\N_0 = \{0,1,2,3,\ldots \}$ 
the set of nonnegative integers, and \Z\ the set of all integers.  
We denote by 
\[
A = \{a_1 < a_2 < \dots < a_k \}
\]
a set $A = \{a_1,a_2,\ldots, a_k\}$ of integers such that $a_1 < a_2 < \dots < a_k$. 
An \emph{interval of integers} is a set of the form 
\[
[a,b] = \{j \in \Z: a \leq j \leq b\} 
\] 
with $a,b \in \Z$ and $a \leq b$.  
The interval $[a,b]$ is \emph{proper} if $a < b$.  
Every proper interval contains at least two integers.  
The \emph{length} of the interval $[a,b]$ is $\ell = b-a$. 
An interval of length $\ell$ contains $\ell+1$ integers. 
An integer $a$ in a subset $A$ of \Z\ is \emph{isolated} 
if it is not a member of a proper interval contained in $A$, that is, 
$a \in A$ but $a-1\notin A$ and $a+1 \notin A$. 

The cardinality of the set $A$ is denoted $|A|$. 
For every nonempty finite or infinite set $X$ of integers, we define 
\[
 \binom{X}{k} = \{A\subseteq X: |A|=k\}. 
\]
We are mostly interested in infinite sets $X$ and, especially, in the sets $\N_0$ and $\Z$.

The $h$-fold sumset of a set $A$ of integers is the set $hA$ of all sums of $h$ 
not necessarily distinct elements of $A$: 
\[
hA = \{a_1+\cdots + a_h: a_i \in A \text{ for all } i \in [1,h]\}.
\]
If $|A|=k$, then $|hA| \leq \binom{h+k-1}{h}$. 
For $n \in \N$, the set $A$ of integers is an \emph{additive $h$-basis for $n$} 
or a \emph{basis of order $h$ for $n$} 
if the sumset $hA$ contains the interval $[0,n]$. 
Bases of order $h$ for $n$ are also called  
\emph{segment bases}  and, in German, \emph{Abschnittsbasen}. 

For every  nonempty finite set $A$ of integers, let $\ell_h(A)$ be the largest integer $n$ 
such that $A$ is a basis of order $h$ for $n$,  
that is, 
\beq            \label{bases:ell}
\ell_h(A) = \max\{ n: [0,n]\subseteq hA\}.
\eeq
If $0 \notin hA$, then $\ell_h(A)$ is undefined. 
If $0 \in A$, then $\ell_1(A) = \max\{ n: [0,n]\subseteq A\}$. 
If $A$ is a finite set of nonnegative integers with $0,1 \in A$ and $a = \max(A)$, 
then it follows from a fundamental theorem of additive number theory 
(Nathanson~\cite{nath72e}) that there are integers
$h_0$ and $d$ such that $\ell_h(A) = ha-d$ for all $h \geq h_0$. 
Here is another example. 
\bt
If $A = \{0\} \cup \{2^i:i \in \N_0\}$,
then
\[
\ell_h(A) = 2^{h+1}-2.
\]
\et 

\begin{proof} 
The proof is by induction on $h$. 
The condition $0 \in A$ implies that if $n$ is a  sum of  $h$  elements $A$, 
then there is a representation of $n$ 
as a  sum of  $h$  elements $A$ with distinct positive summands.  
(For example, $10 = 1+1+4+4 = 0 + 0 + 2 + 8$.) 
It follows that $n \notin hA$ if the binary representation of $n$ has more than $h$ summands.

For $h=1$, we have $[0,2] \subseteq A$ and $3 \notin A$ and so 
$\ell_1(A) = 2=2^2 - 2$. 

For $h=2$, we have $[0,6] \subseteq 2A$ and $7 \notin 2A$ and so 
$\ell_2(A) = 6=2^3 - 2$. 

Let $h \geq 1$.  
If $\ell_{h}(A) = 2^{h+1}-2$, then 
\[
[0, 2^{h+1}-2] \subseteq hA \subseteq (h+1)A
\]
and 
\[
2^{h+1}-1 = \sum_{i=0}^{h} 2^i \in (h+1)A. 
\]
We have 
\[
2^{h+1} + [0, 2^{h+1}-2] = [2^{h+1} , 2^{h+2}-2]  \subseteq (h+1)A  
\]
and so 
\[
 [0, 2^{h+2}-2]  \subseteq (h+1)A. 
\]
However, $2^{h+2}-1 \notin (h+1)A$ because $2^{h+2}-1 = \sum_{i=0}^{h+1}2^i$ 
has $h+2$ positive summands and so $\ell_{h+1}(A) = 2^{h+2}-2$.
This completes the proof. 
\end{proof}

For every nonempty set $X$ of integers and for all positive integers $h$ and $k$, let 
\beq               \label{bases:L}
L_{X,h}(k) = \left\{\ell_h(A): A \in \binom{X}{k}  \right\}.
\eeq

\bt
Let $X$ be a nonempty set of integers and let  $h,k \in \N$.
\benum
\item[(i)]
If $0 \notin X$, then 
\[
L_{X,h}(1)  = L_{X,1}(k) = \emptyset. 
\]
If $0 \in X$, then 
\[
L_{X,h}(1)  = \{0\} 
\]
and
\[
L_{X,1}(k)   = \{n\} 
\]
where $n$ is the largest integer such that $n \leq k-1$ and $[0,n]\subseteq X$. 
\item[(ii)]
For $h=k=2$,  if $0,1 \in X$, then 
\[
L_{X,2}(2)  =  \{0,2\} 
\] 
and if $1 \notin X$ or if $1\in X$ and $0 \notin X$, then 
\[
L_{X,2}(2)   = \{0\}. 
\] 
\eenum
\et

\begin{proof}
The calculation of $L_{X,h}(k)$ is straightforward if $h=1$ or $k=1$.

Let $h=k=2$ and let  $A = \{a_1,a_2\}$ with $a_1 < a_2$.   
Then $2A = \{2a_1,a_1+a_2,2a_2\}$ with 
$2a_1 < a_1+a_2 < 2a_2$.  If $0 \notin 2A$, then $\ell_2(A)$ is undefined.  Suppose $0 \in 2A$.

If $2a_2=0$, then $0 = a_2  \in A$ and $a_1 < a_2 = 0$ and so $\ell_2(A) = 0$. 

If $a_1+a_2=0$, then $a_1 = -a_2 < 0$ and $2a_2 \geq 2$ and so $\ell_2(A) = 0$. 

If $2a_1=0$, then $a_1 = 0$ and $2A = \{0,a_2,2a_2\}$.
If $1 \in A$ and $a_2=1$, then $\ell_2(A) = 2$.  If $1  \notin A$, then 
and $\ell_2(A) = 0$. 
This completes the proof. 
\end{proof}

The problem is to understand the sets  $L_{X,h}(k)$  for $h = 2$ and  all $k \geq 3$ 
and for all $h\geq 3$ and $k \geq 2$. 

\newpage

For the classical case $X = \N_0$, we write 
\[
L_h(k) = L_{\N_0,h}(k). 
\] 

Here is some data:
\begin{align*}
L_h(1) & = \{ 0  \} \hspace{1.05cm}\text{for all $h \in \N$} \\
L_h(2) & = \{ 0, h \} \qquad \text{for all $h \in \N$} \\
L_3(3) & = \{ 0, 3,6,7  \} \\
L_3(4) & = \{ 0,  3, 6, 7, 9, 10, 11, 12, 14, 15 \} \\
L_3(5) & = \{ 0,  3, 6, 7\} \cup [ 9,24] \\
L_3(6) & = \{ 0,  3, 6, 7\} \cup [ 9,36] \\
L_4(3) & = \{ 0, 4,8,10 \} \\ 
L_4(4) & = \{ 0,  4,8,10,12,14,15,16,17,18,20,22, 23, 24, 26 \} \\
L_5(3) & = \{ 0,  5, 10, 13, 14 \} \\
L_5(4) & = \{ 0, 5, 10, 13, 14, 15\} \cup [ 18,  32 ] \cup \{34, 35 \} \\
L_5(5 ) & = [0,44] \setminus \{ 1,2,3,5,6,7,9,11,13,19,42,43 \} \\
L_5(6 ) & =[0,70] \setminus \{ 1,2,3,5,6,7,9,11,13,19\} \\
L_5( 7) & = [0,108] \setminus \{ 1,2,3,5,6,7,9,11,13,19 \} \\
L_6(3) & = \{ 0,  6, 12, 16, 18 \} \\
L_6(4) & = \{ 0,  6, 12, 16, 18, 22, 23, 24, 26, 27, 28, 29 \} \\
L_7(3) & = \{ 0,  7, 14, 19,  22, 23 \} \\
L_7(4) & = \{ 0,  7, 14, 19, 21, 22, 23, 26, 27, 28 \}
\end{align*}

\bt                                                                  \label{bases:theorem:include}
Let $X$ be a nonempty set of integers.  
\benum
\item                    \label{bases:include-1}
If $X$ is   infinite, then 
\[
L_{X,h}(k) \subseteq L_{X,h}(k+1)  
\] 
 for  all positive integers $h$ and $k$.
\item                \label{bases:include-2}
If $\N_0 \subseteq X$, then 
\[
 L_{X,h}(k) \subsetneq L_{X,h}(k+1)  
\] 
 for  all positive integers $h$ and $k$. 
\item                \label{bases:include-3}
There exist nonempty finite and infinite sets $X$ of integers such that 
\[
L_{X,h}(k) = L_{X,h}(k+1) 
\] 
 for  all positive integers $h$ and $k$.
\item                 \label{bases:include-4}
There exist finite sets of integers such that, 
\[
L_{X,h}(k) \setminus L_{X,h}(k+1)  \neq \emptyset 
\] 
 for  all positive integers $h$ and $k$.
\eenum
\et

\begin{proof}
Let $n \in L_{X,h}(k)$.  The set $X$ contains a subset $A$ such that $|A| = k$ and 
$\ell_h(A) = n$, that is, $[0,n] \subseteq hA$ but $n+1 \notin hA$.

If $h = 1$ and $n \in L_{X,1}(k)$, then $[0,n] \subseteq A$ and $n+1 \notin A$.  
Because $X$ is an infinite set, there exists $a' \in X$  such that $|a'| \geq  n+2$ 
and so $A' = A \cup \{a'\}$ is a subset of $X$  with $|A'| = k+1$.  
We have $[0,n] \subseteq hA \subseteq hA'$ and $n+1 \notin A'$ and so 
$n= \ell_1(A') \in L_{X,1}(k+1)$.   Therefore, $L_{X,1}(k) \subseteq L_{X,1}(k+1)$. 

Let $h \geq 2$ and 
\[
v =  \max\{|a|:a\in A\}.
\] 
Then  $hA \subseteq [-hv,hv]$ and so $n \leq hv$. 
Because $X$ is an infinite set, there exists $a' \in X$  such that 
\[
|a'| \geq (h-1)v+n+2 > v
\] 
and so $a' \notin A$.  The set $A' = A \cup \{a'\}$ is a subset of $X$  with $|A'| = k+1$ 
and $[0,n] \subseteq hA \subseteq hA'$.  For all  $x \in hA'\setminus hA$, we have 
$x = b+ra'$ for some $r \in [1,h]$ and 
\[
b \in (h-r)A \subseteq [-(h-1)v,(h-1)v].
\]
It follows that 
\[
|x| \geq r|a'|- |b| \geq |a'|- |b| \geq (h-1)v+n+2 - (h-1)v = n+2
\]
and so $n+1 \notin hA'$ and $\ell_h(A') = n \in  L_{X,h}(k+1)$.  
Therefore, $L_{X,h}(k) \subseteq L_{X,h}(k+1)$.  
This proves~\eqref{bases:include-1}. 

Suppose that $X$ contains $\N_0$. 
There exists $A \in \binom{X}{k}$ such that 
\[
n_0 = \max L_{X,h}(k) = \ell_h(A).
\]  
Then $n_0 \in hA$ but $n_0+1\notin hA$.  
There exist $a_1,\ldots, a_h \in A$ such that 
$a_1 \leq \cdots \leq a_h$ and $n_0 = a_1+\cdots + a_h$.  
Moreover, $n_0+1 \notin hA$ implies $a_h+1 \notin A$.  Also, $n_0 \geq 0$ implies 
$a_h \geq 0$ and so $a_h+1 \in \N_0 \subseteq X$.  
The set $A' = A \cup \{a_h + 1 \} \subseteq X$ satisfies $|A'|=k+1$ and 
$[0,n_0+1] \subseteq hA'$.  It follows that $\ell_h(A') \geq n_0+1 > \max L_{X,h}(k)$ and 
so 
\[
\ell_h(A') \in L_{X,h}(k+1) \setminus L_{X,h}(k)
\] 
This proves~\eqref{bases:include-2}.  

If $X$ is any nonempty set of nonpositive integers that contains 0,
 then $L_{X,h}(k) = \{0\}$ for all 
$h$ and $k$.  This proves~\eqref{bases:include-3}. 

Let $X=[0,k]$.  The unique subset of $X$ of size $k+1$ is $X = [0,k]$.   
We have $hX = [0,hk]$ and $\ell_h(X) = hk$ and so $L_{X,h}(k+1) = \{hk\}$. 
The set $A =  [0,k-1]$ is a subset of $X$ of size $k$ with $hA=[0,h(k-1)]$ 
and 
\[\ell_h(A) = h(k-1) \in L_{X,h}(k) \setminus L_h(k+1).
\] 
This proves~\eqref{bases:include-4}. 
\end{proof}

\bprob[*] 
For positive integers $h$ and $k$, describe the sets $L_h(k)$. 
Similarly, describe the sets $L_{X,h}(k) $ for infinite sets $X$ of integers. \\
Discover the beautiful patterns hidden in these sets. 
\eprob

\bprob
For fixed $h\in \N$, compute the sequence $\left( L_{h}(k)  \right)_{k=1}^{\infty}$. 
AI should be able to compute the finite sequence $\left( L_{h}(k)\right)_{k=1}^K$ for large $K$ 
and then predict (find a formula?) for sufficiently large  values of $L_{h}(k)$. 
To understand $L_{h}(k)$ for small values of $k$ might be more difficult. 
\eprob

\bprob
For positive integers $h$ and $k$, find necessary and sufficient conditions 
on an infinite set $X$ of integers such that 
\[
L_{X,h}(k) \subsetneq L_{X,h}(k+1). 
\]
\eprob

\bprob
Find necessary and sufficient conditions 
on an infinite set $X$ of integers such that 
\[
L_{X,h}(k) \subsetneq L_{X,h}(k+1) 
\]
for all $h$ and $k$ or for all sufficiently large $h$ and $k$. 
\eprob

\bprob
For positive integers $h$ and $k$, find necessary and sufficient conditions 
on an infinite set $X$ of nonnegative integers such that 
\[
L_{X,h}(k) \subsetneq L_{X,h}(k+1). 
\]
\eprob

\bprob
Find necessary and sufficient conditions 
on an infinite set $X$ of nonnegative  integers such that 
\[
L_{X,h}(k) \subsetneq L_{X,h}(k+1) 
\]
for all $h$ and $k$ or for all sufficiently large $h$ and $k$. 
\eprob

\bprob[*]
Consider finite sets that can contain negative numbers. 
We have $L_{h}(k) \subseteq L_{\Z,h}(k)$.  
Describe the set $L_{\Z,h}(k) \setminus L_{h}(k)$.
\eprob

An infinite set  $X$ contains infinitely many $k$-element subsets,  
but there is an immediate reduction of the  computation of   $L_{X,h}(k)$ 
to an examination of only finitely  many $k$-element subsets of $X$.

For $k \in \N$, let  $\mca_h(k)$ be the set of all sets 
$A = \{a_1, a_2, \ldots, a_k\} \in \binom{\mathbf{N}_0}{k}$ such that 
\[
0 = a_1 < a_2 < \cdots < a_k 
\]
and
\[
a_j+1 \leq  a_{j+1}  \leq ha_j+1 \qquad \text{ for all $j \in [1,k-1]$}.
\]
Thus,  $\mca_h(1) = \{\{0\}\}$, $\mca_h(2) =  \{\{0,1 \}\}$, and 
\[
\mca_h(3) =  \{\{0,1,a_3 : a_3 \in [2,h+1] \}\}. 
\] 

\bt              \label{bases:theorem:compress}
For all positive integers $h$ and $k$, 
\[
L_{h}(k)=  \left\{\ell_h(A): A \in \mca_h(k) \right\} \cup  \bigcup_{j=1}^{k-1} L_h(j) .
\]
\et

\begin{proof}
Let $A = \{a_1, a_2, \ldots, a_k\}$ be a set of nonnegative integers with 
\[
0 = a_1 < a_2 < \cdots < a_k.
\] 
If 
\[
a_{i+1} \leq ha_i+ 1   
\]
for all $i \in [1,k-1]$, then $A \in \mca_h(k)$ and  $\ell_h(A) \in \left\{\ell_h(A): A \in \mca_h(k) \right\}$. 

Suppose that 
\[
a_{i+1} \geq ha_i+ 2 
\]
for some $i \in [1,k-1]$.  
Let $j$ be the smallest integer in $[1,k-1]$ such that $a_{j+1} \geq ha_j+ 2$ 
and let $A' = \{a_1, a_2, \ldots, a_j\} \in \binom{\mathbf{N}_0}{j}$. 
Then $A \setminus A' = \{a_{j+1},a_{j+2},\ldots, a_k\}$.  
We have $hA' \subseteq hA$ and $x \leq ha_j$ for all $x \in hA'$.  
If $x \in hA\setminus hA'$, then there exists $r \in [1,h]$ 
such that $x = u+v$ with $u \in (h-r)A'$ and $v \in r(A\setminus A')$.
Then 
\[
x  \geq v \geq a_{j+1} \geq ha_j+2 
\]
and so $ha_j+1 \notin hA$. Therefore, $\ell_h(A) \leq ha_j$.  
It follows that $x \in hA \cap [0,ha_j]$ if and only if $x\in hA'$ 
and so 
\[
\ell_h(A) = \ell_h(A') \in L_h(j) \subseteq  \bigcup_{j=1}^{k-1} L_h(j).
\]
This completes the proof. 
\end{proof}


\section{Interval bases of order $h$ for $n$}        \label{bases:section:four} 

The set $A$ of integers is an \emph{interval $h$-basis for $n$} 
or an \emph{interval basis of order $h$ for $n$} if the sumset $hA$ 
contains an interval $[c, c+n]$ for some $c \in \Z$. 
For every  nonempty finite set $A$ of integers, let $\ell^{\sharp}_h(A)$ be the largest integer $n$ 
such that $A$ is an interval $h$-basis for $n$,  
that is,  the largest integer $n$ such that $[c,c+n] \subseteq hA$ for some $c \in \Z$.  

The function  $\ell^{\sharp}_h(A)$ is invariant under reflection and translation.  
For $b \in \Z$ and  $\varepsilon = -1$, we define the \emph{reflection} 
\[
-A  = \{ -a:a \in A\}
\]
and the \emph{translation}  
\[
 A + b=  \{a+b:a \in A\}.
\]
If  $\varepsilon = \pm1$ and 
\[
B =   \varepsilon A  +b= \{ \varepsilon a +b:a \in A\}
\] 
then
\[
hB = h(  \varepsilon A +b) =    \varepsilon (hA) +hb
\]
and so 
\beq          \label{bases:affine}
\ell^{\sharp}_h(A) = \ell^{\sharp}_h(  \varepsilon A +b).
\eeq 
For example, if
\[
A = \{0,1,3\} 
\]
and 
\[
B = 2 -A  = \{-1,1,2\} 
\]
then  
\[
2A = \{0,1,2,3,4,6\}
\]
and 
\[
2B = \{-2,0,1,2,3,4\} = 4-2A. 
\]
We have  $\ell^{\sharp}_h(A) = \ell^{\sharp}_h(b + \varepsilon A) =2$. 

For any nonempty set $X$ of integers and all positive integers $h$ and $k$, let 
\[
L^{\sharp}_{X,h}(k) = \left\{\ell^{\sharp}_h(A): A \in \binom{X}{k}  \right\}.
\] 
For sets of nonnegative integers, we write 
\[
L^{\sharp}_h(k)= \left\{\ell^{\sharp}_h(A): A \in \binom{\N_0}{k}  \right\}.
\]

\bprob[*] 
For positive integers $h$ and $k$, describe the sets $L^{\sharp}_h(k)$.  
Similarly, describe the sets $L^{\sharp}_{X,h}(k)$ for infinite sets $X$ of integers. \\
What patterns are hidden in these sets?  
\eprob

\bprob[*]
What  are the relations between the sets  $L^{\sharp}_h(k)$ and $L_{\Z,h}(k)$? 
\eprob

For positive integers $h$ and $k$, we consider the following four 
additive number theoretic functions defined by $h$-fold sumsets of sets of size $k$.

\benum
\item
Let $n^{\flat}_h(k)$ be the largest integer $n$ such that 
there exists a set $A$ of $k$ nonnegative integers with $[0,n] \subseteq hA$: 
\[
n^{\flat}_h(k) = \max L_{h}(k) = \max\left\{ \ell_h(A): A \in  \binom{\N_0}{k} \right\}.
\] 

\item
Let $n^{\sharp}_h(k)$ be the largest integer $n$ such that 
there exists a set $A$ of $k$ nonnegative integers with $[c, c+n]\subseteq hA$ for some integer $c$: 
\[
n^{\sharp}_h(k) = \max\left\{ \ell^{\sharp}_h(A): A \in  \binom{\N_0}{k} \right\}.
\]

\item
Let $m^{\flat}_h(k)$ be the largest integer $n$ such that 
there exists a set $A$ of $k$  integers such that $[0,n] \subseteq hA$:  
\[
m^{\flat}_h(k) = \max\left\{ \ell_h(A): A \in  \binom{\Z}{k} \right\}.
\] 

\item
Let $m^{\sharp}_h(k)$ be the largest integer $n$ such that there exists a set $A$ of $k$  integers 
with $[c, c+n]\subseteq hA$ for some integer $c$:  
\[
m^{\sharp}_h(k) = \max\left\{ \ell^{\sharp}_h(A): A \in  \binom{\Z}{k} \right\}.
\] 
\eenum

For all $k \in \N$, we have 
\[
n^{\flat}_1(k) = n^{\sharp}_1(k) = m^{\flat}_1(k) = m^{\sharp}_1(k) = k - 1 
\]
and the interval  $[0,k-1]$ satisfies this maximal condition.   
For all $h \in \N$, we have 
\[
n^{\flat}_h(1) = n^{\sharp}_h(1) =  m^{\flat}_h(1) = m^{\sharp}_h(1) = 1   
\]
and the set $\{0\}$ satisfies this maximal condition.   
We investigate the  functions $n^{\flat}_h(k)$,  $n^{\sharp}_h(k)$, 
 $m^{\flat}_h(k)$, and $m^{\sharp}_h(k)$ for $h \geq 2$ and $k \geq 2$.

\bt              \label{bases:theorem:strict} 
For fixed $h \in \N$, the  sequences $(n^{\flat}_h(k))_{k=1}^{\infty}$ , 
$(m^{\flat}_h(k))_{k=1}^{\infty}$, 
 $(n^{\sharp}_h(k))_{k=1}^{\infty}$, and $(m^{\sharp}_h(k))_{k=1}^{\infty}$
are strictly increasing.  
\et

\begin{proof}
Fix $h\in \N$.  For $k \in \N$, choose $A \in \binom{\Z}{k}$ and $c \in \N_0$ 
such that  $\ell^{\sharp}_h(A) = m^{\sharp}_h(k) = n$ 
and $[c,c+n] \subseteq hA$.  
There exist integers $a_1,\ldots, a_{h-1}, a_h \in A$ with $a_1 \leq \cdots \leq a_{h-1} \leq a_h$ 
such that 
\[
a_1+\cdots + a_{h-1} + a_h = c+ n. 
\]
The maximality of $n$ implies  
\[
a_1+\cdots + a_{h-1} + (a_h+1) = c+ n+1 \notin hA 
\]
and so  $a_h+1 \notin A$.  
It follows that the set  $A' = A \cup \{a_h+1\}$ satisfies $|A'| = k+1$ amd $c+n+1 \in hA'$ and so 
\[
[c,c+n+1] = [c,c+n] \cup \{c+n+1\} \subseteq hA \cup \{c+n+1\} \subseteq hA'. 
\] 
Therefore, $m^{\sharp}_h(k+1) \geq \ell^{\sharp}_h(A') \geq n+1 > m^{\sharp}_h(k) $. 

This argument with $c=0$ proves $m^{\flat}_h(k+1) > m^{\flat}_h(k) $ and the argument 
also works with the functions $n^{\sharp}_h(k)$ and $n^{\flat}_h(k)$.  
This completes the proof. 
\end{proof}

\bprob[*]
Compute the sequences  $(n^{\flat}_h(k))_{k=1}^{\infty}$ , $(m^{\flat}_h(k))_{k=1}^{\infty}$, 
 $(n^{\sharp}_h(k))_{k=1}^{\infty}$, and $(m^{\sharp}_h(k))_{k=1}^{\infty}$. 
Construct an algorithm to decide, for every positive integer $w$,
 if there exist $h$ and $k$ such that $w = n^{\flat}_h(k)$ or $w = m^{\flat}_h(k)$ 
 or $w = n^{\sharp}_h(k)$ or $w = m^{\sharp}_h(k)$.  
\eprob

\bprob[*]
For $h, d \in \N$, describe the difference sequences 
$(n^{\flat}_h(k+d) - n^{\flat}_h(k))_{k=1}^{\infty}$, $(m^{\flat}_h(k+d) - m^{\flat}_h(k))_{k=1}^{\infty}$, 
$n^{\sharp}_h(k+d) - (n^{\sharp}_h(k))_{k=1}^{\infty}$, and 
$(m^{\sharp}_h(k+d) - m^{\sharp}_h(k))_{k=1}^{\infty}$.   
Are these sequences bounded for some or all $d$? 
\eprob

\bprob 
It is natural to ask how the use of negative numbers changes the size of the maximal 
interval $[0,n]$ contained in an $h$-fold sumset. 
For integers $h \geq 2$ and $k \geq 2$, are the following statements true or false?
\benum 
\item
If $A \in \binom{\Z}{k}$ with $\min(A) < 0$, then 
\[
\ell_h(A) \leq n^{\flat}_h(k).
\]

\item
If $A \in \binom{\Z}{k}$ with $\min(A) < 0$, then 
\[
\ell_h(A) < n^{\flat}_h(k).
\]
\eenum 

\eprob

\section{A generalization of Sidon sets}         \label{bases:Sidon}
The study of additive $h$-bases for $n$ has consisted mostly of 
(still unsuccessful) attempts to compute 
or to estimate asymptotically the  counting function  $n^{\flat}_h(k)$, 
often by constructing explicit families of such sets.   
Also investigated is the number of subsets   of $[0,n]$ that are $h$-bases for $n$ 
(for calculations with $h=2$, see OEIS A066062) and the size, denoted $k_h(n)$, 
of the smallest subset of $[0,n]$ 
that is an $h$-basis for $n$ (for calculations of $k_2(n)$, see OEIS A066063).  The other 
counting functions constructed in Section~\ref{bases:section:four} are largely unstudied.  

Little is known about the possible arithmetical structures of interval $h$-bases for $n$.
In Section~\ref{bases:section:structure} we construct a new class of interval 
2-bases. The proof uses a generalization of the concept of Sidon set.

A \emph{Sidon set} is a finite or infinite set $A = \{a_i: i \in I\}$ of positive integers 
such that, for all ${r_1},{r_2},{s_1},{s_2} \in I$, we have 
\[
a_{r_1}+a_{r_2} = a_{s_1} + a_{s_2}
\] 
if and only if 
\[
\{{r_1},{r_2}\} = \{ {s_1},{s_2}\}.
\]
Equivalently, if 
\beq                    \label{bases:Sidon-1}
\{{r_1},{r_2}\} \neq \{ {s_1},{s_2}\} 
\eeq
then 
\[
\left| (a_{r_1}+a_{r_2})-(a_{s_1}+a_{s_2}) \right| \geq 1. 
\]
Let $\Delta > 0$.  The set $A$ has a $\Delta$-separated sumset if~\eqref{bases:Sidon-1} impllies 
\[
\left| (a_{r_1}+a_{r_2})-(a_{s_1}+a_{s_2}) \right|  \geq \Delta. 
\]
More generally, for $h \geq 2$, a \emph{$B_h$-set} is a set $A = \{a_i:i \in I\}$ of integers such that, 
for all $r_1,\ldots, r_h, s_1,\ldots, s_h \in I$, conditions 
\beq                    \label{bases:Sidon-2a}
r_1 \leq \cdots \leq  r_h \qqand  s_1\leq \cdots \leq s_h
\eeq
and 
\beq                    \label{bases:Sidon-2b}
r_i \neq s_i \qquad \text{for some $i \in [1,h]$}
\eeq 
imply 
\[
\left| \sum_{i=1}^h a_{r_i} - \sum_{i=1}^h a_{s_i} \right| \geq 1.
\]
The set $A$ has a \emph{$\Delta$-separated $h$-fold sumset} if 
conditions~\eqref{bases:Sidon-2a} and~\eqref{bases:Sidon-2b} imply 
\[
\left| \sum_{i=1}^h a_{r_i} - \sum_{i=1}^h a_{s_i} \right| \geq \Delta.
\]
Thus, a $B_h$-set is a set with a $1$-separated $h$-fold sumset. 
If the set $A$ has a \emph{$\Delta$-separated $h$-fold sumset}, then 
$A$ has a \emph{$\Delta$-separated $h'$-fold sumset} for all $h' \in [1,h]$.

Let $H \geq 2$.  The set $A$ is \emph{ $\Delta$-separated of level $H$} if 
\benum
\item
$A$ has a $\Delta$-separated $h$-fold sumset for all $h \in [1,H]$, and 
\item
if $h_1,h_2 \in [1,H]$ with $h_1 \neq h_2$, then 
\[
|u_1-u_2| \geq \Delta \qquad \text{for all $u_1 \in h_1A$ and $u_2 \in h_2A$.}
\]
\eenum
 The set $A$ is \emph{ $\Delta$-separated of level $\infty$} if $A$ is 
 $\Delta$-separated of level $H$ for all $H \in \N$.  
 
 \bl                         \label{bases:lemma:separation}
 The set $A = \{a_i:i\in I\}$ is $\Delta$-separated of level $\infty$ if and only if, 
 for all finite subsets $I_1$ and $I_2$  of $I$ with $I_1 \neq I_2$, we have 
\beq                     \label{bases:separation}
 \left| \sum_{i  \in I_1}a_{r_i} - \sum_{i \in I_2} a_{s_i} \right| \geq \Delta. 
 \eeq
 \el 

\begin{proof}
This follows immediately from the definition of $\Delta$-separated of level $\infty$. 
\end{proof}

For the construction of the special class of interval $2$-bases for $n$ in Theorem~\ref{bases:theorem:construction}, 
we need only the existence 
of arbitrarily large finite  $\Delta$-separated sets of level $3$ for all $\Delta \in \N$.   
Every subset of a  $\Delta$-separated set of level $H$ is a  $\Delta$-separated of level $H$, 
and so the following simple result suffices. 

\bt               \label{bases:theorem:g-Delta}  
For $g \geq 2$ and $g^{i_0} \geq \Delta$, the set  $\{g^i:i\geq i_0\}$ is  
$\Delta$-separated of level $\infty$. 
\et

\begin{proof}
Let $I_1$ and $I_2$ be finite subsets of the set $\{i\in \N:i\geq i_0\}$  with $I_1 \neq I_2$ 
and let  
\[
u_1 = \sum_{i_1 \in I_1} g^{i_1} = g^{i_0} \sum_{i_1 \in I_1} g^{i_1-i_0}
\]
and 
\[
u_2 = \sum_{i_2 \in I_2} g^{i_2} = g^{i_0} \sum_{i_2 \in I_2} g^{i_2-i_0}.
\]
The uniqueness of the $g$-adic representation implies  $u_1 \neq u_2$ and so 
\[
|u_1-u_2| =  g^{i_0} \left| \sum_{i_1 \in I_1} g^{i_1-i_0} -  \sum_{i_2 \in I_2} g^{i_2-i_0} \right| 
\geq g^{i_0} \geq \Delta.
\]
This completes the proof.  
\end{proof}

\bprob[*]
Determine a reasonable notion of ``density'' and construct dense finite and infinite 
sets of positive integers that are $\Delta$-separated of level $H$.
\eprob

\section{A construction of interval $2$-bases for $n$} \label{bases:section:structure} 
Recall that, for every  nonempty finite set $A$ of integers,  $\ell^{\sharp}_h(A)$ is 
 the largest integer $n$ such that $[c,c+n] \subseteq hA$ for some $c \in \Z$.

\bt             \label{bases:theorem:construction}
For all integers $c$, $n$, and $d$  with $n \geq 1$ and $d \geq 2$, 
there is a set $A$ of integers with $|A| = 2k+2$ such that 
\benum
\item[(i)]
$\ell^{\sharp}_2(A) = n$ and $[c,c+n]$ is the unique interval of length $n$  
contained in the sumset $2A$ 
\item[(ii)] 
 $|x-y| \geq d$ for all $x,y \in 2A$ with 
$x \neq y$ and $x \notin [c,c+n]$ . 
\eenum 
\et

Thus, every element of the set $2A\setminus [c,c+n]$ is isolated.

\begin{proof}
Let 
\[
\Delta  =   |c|+d+ 2n.
\]
Let $(a_j)_{j=0}^n$ be a strictly increasing sequence of positive integers such that  
\beq              \label{bases:Delta-0}
a_0 \geq \Delta 
\eeq
and 
\beq              \label{bases:Delta}
a_j > 2a_{j-1} + |c| + n 
\eeq 
for all $j \in [1,n]$.  Consider the sets 
\begin{align*}
A _0 & =  \{ c+j + a_j:j \in [0,n] \} \\ 
A_1  & =  \{-a_j :j \in [0,n] \} 
\end{align*}
and 
\[
A = A_0 \cup A_1. 
\] 
Then $|A_0| = |A_1| = n+1$.  The inequality 
\[
c+j +a_j \geq c+a_0 \geq d + 2n > 0 > -a_j
\]
implies $A_0 \cap A_1 = \emptyset$ and $|A| = 2n+2$.  

We have the sumset 
\[
2A = B_{0,1} \cup B_{0,0} \cup B_{1,1}
\]
where 
\begin{align*}
B_{0,0}  &  = 2A_0   = \{2c + i+j +  a_i + a_j   : i,j \in [0,n]\} \\ 
B_{0,1} &  = A_0 + A_1   = \{c + j + a_j - a_i : i,j \in [0,n]\} \\
B_{1,1} & = 2A_1   = \{  -a_i - a_j : i,j \in [0,n]\}.
\end{align*} 

For all $i,j \in [0,n]$, let 
\[
b_{i,j} = c+j + a_j -a_i \in B_{0,1}.
\] 
Then 
\[
b_{j,j} = c+j + a_j -a_j = c+j  \in B_{0,1} \subseteq 2A
\]
for all $j \in [0,n]$ and so $[c,c+n]   \subseteq  2A$.  
Thus, $\ell^{\sharp}_h(A) \geq n$.

To obtain the element isolation inequality~(ii), we choose a set 
$(a_j)_{j=0}^n$ that is $\Delta$-separated of level $3$ 
and satisfies inequalities~\eqref{bases:Delta-0} and~\eqref{bases:Delta}.  
For example, applying Theorem~\ref{bases:theorem:g-Delta} 
with 
\[
g=\Delta   \qqand a_j = g^{j+1}  
\]  
for all $j \in [0,n] $, we have 
\[
a_0 = g  = \Delta \geq 4 
\]
 and 
\[
a_j = g^{j+1} > 3g^j = 2g^{j} +g^{j} \geq  2a_{j-1} + \Delta >  2a_{j-1} + 2|c| + d + 2n  
\]
for all $j \in [1,n]$. 

In each of the following six cases, we apply inequality~\eqref{bases:separation} of Lemma~\ref{bases:lemma:separation}. 
\benum 
\item
If $x,y \in B_{0,0}$ and $x \neq y$, then there exist $\{r_1,r_2, s_1,s_2\} \subseteq [0,n] $
with $\{r_1,r_2\} \neq \{s_1,s_2\} $ such that 
\[
x = 2c+ r_1+ r_2 + a_{r_1} + a_{r_2} 
\]
\[
y = 2c+s_1+ s_2 + a_{s_1} + a_{s_2}.
\]
Then 
\[
x-y = r_1+ r_2 - s_1-  s_2  + a_{r_1} + a_{r_2} -  a_{s_1}  - a_{s_2}  
\]
and 
\begin{align*}
|x-y| 
& = \left|r_1+ r_2 - s_1-  s_2  + a_{r_1} + a_{r_2} -  a_{s_1}  - a_{s_2} \right| \\
& \geq \left| \left(  a_{r_1} + a_{r_2} \right) - \left(a_{s_2} + a_{s_2}  \right) \right| 
- \left| r_1+ r_2 -  s_1- s_2 )  \right| \\
&  \geq \Delta - 2n \\
& \geq d.
\end{align*}

\item
If $x,y \in B_{1,1}$ and $x \neq y$, then there exist $\{r_1,r_2,s_1,s_2\} \subseteq [0,n]$ 
with $\{r_1,r_2\} \neq \{s_1,s_2\} $ such that 
\[
x= -  a_{s_1} -  a_{s_2}  
\] 
\[
y = -a_{r_1} - a_{r_2}. 
\]
Then  
\[
x-y  =  a_{r_1} + a_{r_2}  - a_{s_1} - a_{s_2}  
\]
and 
\begin{align*}
|x-y| & = \left| a_{r_1} + a_{r_2}  - a_{s_1} - a_{s_2}   \right| \\
&  \geq \Delta  \\
& \geq d.
\end{align*}

\item
If $x \in B_{0,0}$ and $y \in B_{1,1}$, then $y < 0 < x$ 
and there exist $\{r_1,r_2, r_3,r_4\} \subseteq [0,n]$ such that 
\begin{align*}
x & = 2c+r_1+ r_2 + a_{r_1} + a_{r_2} \\
y & = -a_{r_3} - a_{r_4}. 
\end{align*}
Recalling inequality~\eqref{bases:Delta-0}, we obtain  
\begin{align*}
|x-y| & =  2c+r_1+ r_2 + a_{r_1} + a_{r_2} +  a_{r_3} + a_{r_4} \\
& \geq 4a_0 - 2|c|  \\
& \geq 4\Delta - 2|c|  \\
& \geq  d.
\end{align*}

\item
If  $x \in B_{0,1}$ and $y \in B_{0,0}$, 
then there exist $\{r_1,r_2, r_3,s_1\} \subseteq [0,n]$ such that 
\begin{align*} 
x & =  c+r_1+  a_{r_1}  - a_{s_1} \\
y & =  2c+s_2+ s_3 + a_{s_2} + a_{s_3} 
\end{align*} 
and 
\[
x-y = -c +r_1-s_2-s_3 + a_{r_1} - \left(  a_{s_1} + a_{s_2} + a_{s_3}  \right). 
\]
If $ \{ r_1\} = \{s_1,s_2, s_3\}$, then $r_1 = s_1 = s_2 = s_3 $ and 
\[
x-y =  -c  - r_1 - 2 a_{r_1} 
\]
and so 
\begin{align*} 
|x-y|  & = | -c  - r_1 - 2 a_{r_1} | \\
& \geq 2a_{r_1} - |c|- r_1 \\
& \geq 2a_0 - |c| - n \\
& \geq  2\Delta - |c|-n \\
& \geq d.
\end{align*}
If $ \{ r_1\} \neq \{s_1,s_2, s_3\}$, then
\begin{align*}
|x-y| 
& = \left|  -c +r_1-s_2-s_3 + a_{r_1} - \left(  a_{s_1} + a_{s_2} + a_{s_3}  \right) \right| \\
& \geq \left|  a_{r_1} - \left(  a_{s_1} + a_{s_2} + a_{s_3}  \right)   \right|  -  \left|  -c +r_1-s_2-s_3 \right| \\
&  \geq \Delta - |c| - 2n \\
& \geq d. 
\end{align*} 

\item
If  $x \in B_{0,1}$ and $y \in B_{1,1}$, 
then there exist $\{r_1,r_2, r_3,s_1\} \subseteq [0,n]$ such that 
\begin{align*} 
x & =  c+r_1+  a_{r_1}  - a_{s_1} \\
y & =  -a_{r_2} - a_{r_3} 
\end{align*} 
and 
\[
x-y = c+r_1+  a_{r_1} +  a_{r_2} + a_{r_3}  - a_{s_1}
 \]
If $\{r_1,r_2,r_3\} = \{ s_1\}$, then $r_1 = r_2 = r_3 = s_1$ and 
\[
x-y = c+r_1 -2a_{r_1} 
\]
and   
\begin{align*} 
|x-y|  & = |c+r_1 -2a_{r_1} | \\
& \geq 2a_{r_1} - |c|-n \\
&   > \Delta - |c|-n \\
& \geq d.
\end{align*}
If $\{r_1,r_2,r_3\} \neq \{ s_1\}$, then
\begin{align*}
|x-y| 
& = \left|  c+r_1+  a_{r_1} +  a_{r_2} + a_{r_3}  - a_{s_1} \right| \\
& \geq \left| \left(  a_{r_2} + a_{r_2} + a_{r_3}  \right) - a_{s_1}  
 \right| -|c+r_1 |  \\
&  \geq \Delta -  |c| - n \\
& \geq d.
\end{align*} 

\item 
If $x,y \in B_{0,1}$ and $x \neq y$ and $x \notin [c,c+n]$, 
then there exist  $\{r_1,r_2, s_1, s_2\} \subseteq [0,n]$ with $ r_1 \neq s_1$ such that 
\begin{align*} 
x & = c+r_1+a_{r_1}-a_{s_1} \\
y & = c+s_2+a_{s_2}-a_{r_2}.
\end{align*} 
and 
\[
x-y = r_1 - s_2 + a_{r_1} + a_{r_2}  - a_{s_1}  - a_{s_2}. 
\]
If $\{r_1,r_2\} =\{s_1,s_2\} $, then the inequality $r_1 \neq s_1$ implies $r_1 = s_2$ and 
$r_2=s_1$, and so $x=y$, which is absurd.  Therefore, 
 $\{r_1,r_2\}  \neq \{s_1,s_2\} $ and 
\begin{align*}
|x-y| & = | r_1 - s_2 + a_{r_1} + a_{r_2}  - a_{s_1}  - a_{s_2}| \\ 
& \geq  \left| \left( a_{r_1} + a_{r_2} \right)  - \left( a_{s_1}  + a_{s_2} \right)  \right|  
- \left| r_1 - s_2 \right| \\ 
& \geq \Delta - n \\ 
& \geq d. 
\end{align*} 
\eenum
This completes the proof. 
\end{proof}

Here is a useful observation about translations of  intervals in sumsets. 

\bt           \label{bases:theorem:r-translate}
 Let $A$ be a set of integers. 
\benum
\item[(a)]
  If $[c,c+n] \subseteq hA$, then there exist $r \in [0,h-1]$ 
 and $c_0 \in \Z$ such that $[r,r+n] \subseteq h(A-c_0)$.
 \item[(b)]
For all integers $h \geq 2$ and $k \geq 2$, 
\[
m^{\sharp}_h(k) =  \max\left\{ 
n \in \N: [r,r+n] \subseteq hA \text{ for some } A \in  \binom{\Z}{k} \text{ and } r \in [0,h-1]\right\}. 
\] 
 \eenum
\et

\begin{proof} 
Let $A \subseteq \Z$ with $[c,c+n] \subseteq hA$ for some $c \in \Z$.  
By the division algorithm, we have 
\[
c = hq+r
\]
for some $q \in \Z$ and $r \in [0, h-1]$. 
By~\eqref{bases:affine}, the translate $ A-q  \subseteq \Z$ satisfies  $|A- q| = |A|$ and  
\begin{align*} 
 h(A-q) & = hA - hq \\
& \supseteq [c,c+n] - hq =  [c-hq, c-hq +n] \\
&   = [r,r+n]. 
\end{align*}
Therefore, 
\begin{align*} 
m^{\sharp}_h(k) 
& =  \max\left\{ 
n \in \N: [c,c+n] \subseteq hA \text{ for some } A \in  \binom{\Z}{k} \text{ and } c \in \Z \right\} \\
& =  \max\left\{ 
n \in \N: [r,r+n] \subseteq hA \text{ for some } A \in  \binom{\Z}{k} \text{ and } r \in [0,h-1]\right\}.
\end{align*}
This completes the proof. 
\end{proof}

\bprob
Let  $h \geq 2$,  $q \geq 2$ and  $n \geq 1$ be   integers  and 
let $C = \{c_i:i\in [1,q] \}$ be a set of integers with $c_i - c_{i-1} \geq n+2$ 
for all $i \in [2,q]$.  The intervals $[c_i,c_i+n]$ are pairwise disjoint for $i \in [1,k]$. 
Let $h \geq 2$.  Is it true that  there exists a finite set $A$ such that 
\[
[c, c+n] \subseteq hA 
\]
if and only if $c = c_i \in C$. 
\eprob

\bprob
Let  $h \geq 2$,  $q \geq 2$, $n \geq 1$, and $N_q \geq n+2$  be integers   
and let $C = \{c_i:i\in [1,q] \}$ be a set of integers with $c_i - c_{i-1} \geq  N_q$ 
for all $i \in [2,q]$.  
Is it true that there exists a finite set $A$ such that 
\[
[c, c+n] \subseteq hA 
\qquad
\text{if and only if} \qquad c = c_i \in C 
\]
and 
\[
[c',c'+1]  \subseteq hA \qquad\text{ if and only if }
\]
 \[
[c', c'+1] \subseteq [c_i, c_i+n]\subseteq hA \quad \text{for some $c_i \in C$.} 
\]
If true, compute or estimate the size of $A$ as a function of $h,q,n$, and $N_q$. 
\eprob

\section{Elementary inequalities} 
The relation  $ \binom{\N_0}{k}\subseteq  \binom{\Z}{k}$ 
implies $n^{\sharp}_h(k) \leq m^{\sharp}_h(k)$.
The following lemma removes the inequality.  

\bt 
For all positive integers $h$ and $k$, 
\[
L^{\sharp}_h(k) = L^{\sharp}_{X,h}(k) 
\]
and 
\[
n^{\sharp}_h(k) = m^{\sharp}_h(k).
\]
\et

\begin{proof}
Let $A  \subseteq \Z$.  For all   $b, c \in \Z$ and $n \in \N$, we have 
\[
[c,c+n] \subseteq hA 
\]
if and only if 
\begin{align*} 
[c+hb,c+hb+n] & = [c,c+n]+hb \\
& \subseteq hA + hb  = h(A+b)  
\end{align*} 
and 
\[
 c+n+1 \notin hA 
 \]
 if and only if 
\[
c+hb+n+1 \notin hA+hb = h(A+b)
\]
Thus, the set of lengths of intervals contained in $A$ is the same as the set of lengths
of intervals contained in the translate $A+b$. 

If $A \in  \binom{\Z}{k}$ and $A \notin  \binom{\N_0}{k}$, then $\min(A) = -a_0 < 0$ and 
$A+a_0 \in  \binom{\N_0}{k}$.   
It follows that 
 $L^{\sharp}_{\N_0,h}(k) = L^{\sharp}_{\Z,h}(k)$ and 
\[
n^{\sharp}_h(k) = \max L^{\sharp}_{\N_0,h}(k) 
= \max L^{\sharp}_{\Z,h}(k) = m^{\sharp}_h(k). 
\]
This completes the proof. 
\end{proof}

For integers $h \geq 2$ and $k \geq 2$,  the elementary inequalities 
\begin{align*}
n^{\flat}_h(k) & \leq m^{\flat}_h(k)\\
n^{\flat}_h(k)  & \leq n^{\sharp}_h(k) \\
m^{\flat}_h(k) & \leq  m^{\sharp}_h(k).
\end{align*}
suggest the following questions.

\bprob 
For  $h \geq 2$, describe the set  
\[
\{ m^{\flat}_h(k) - n ^{\flat}_h(k) : k \in \N\}
\]
\eprob

\bprob 
For   $h \geq 2$, describe the set  
\[
\{ n^{\sharp}_h(k) - n ^{\flat}_h(k) : k \in \N\} 
\] 
\eprob 

\bprob 
For  $h \geq 2$, describe the set  
\[
\{ m^{\sharp}_h(k) - m^{\flat}_h(k): k \in \N\}. 
\]
\eprob

\section{The additive diameters $d_h(k)$ and  $d^{\sharp}_h(k)$}
 The \emph{diameter} of a bounded set $A$ of numbers is 
 \[
 \diam(A) = \sup\{|a-a'| : a,a' \in A\}.
 \] 
 
 \bt
 If $A \in \binom{\N_0}{k}$ and $\ell_h(A) = n^{\flat}_h(k)$, then $\diam(A) \leq n^{\flat}_h(k)$. 

 \et
 
\begin{proof}
Let $A \in \binom{\N_0}{k}$.  If $\ell_h(A) = n^{\flat}_h(k)$, then 
 $\min(A) = 0$ and $\max(A) = \diam(A) = a_k$.   
Because $A \subseteq \N_0$, if $a_k > n^{\flat}_h(k)$, then the representation 
of an integer in the interval $[0,n^{\flat}_h(k)]$ as a sum of $h$ elements of $A$ does not include 
the number $a_k$ as a summand, and so $A' = A\setminus \{a_k\}$ 
satisfies $A' \in \binom{\N_0}{k-1}$ 
and $[0,n^{\flat}_h(k)] \subseteq hA'$.  Therefore, 
\[
n^{\flat}_h(k) \leq \ell_h(A') \leq  n^{\flat}_h(k-1)
\]
which contradicts Theorem~\ref{bases:theorem:strict}.  
It follows that  $\diam(A) = a_k \leq n^{\flat}_h(k)$. 
This completes the proof. 
\end{proof}

\bprob 
Can sets of bounded size and  large diameter   contain large intervals?   
Describe the following sets:  
\[
\left\{ \diam(A): A \in \binom{\\N_0}{k} \text{ and } \ell_h(A) = n^{\sharp}_h(k) \right\} 
\]
and 
\[
\left\{ \diam(A): A \in \binom{\Z}{k} \text{ and } \ell_h(A) = m^{\flat}_h(k) \right\}.
\]
\eprob

\section{Notes}
We have looked at  $n^{\flat}_h(k)$ as a function of $k$ for fixed $h$, but there are also 
 problems and results when considering  $n^{\flat}_h(k)$ as a function of $h$ for fixed $k$,
and also symmetry questions relating $n^{\flat}_h(k)$ and $n^{\flat}_k(h)$.  

There is also an important function dual to $n^{\flat}_h(k)$.   
For positive integers $h$ and $n$, let $k_h(n)$ be the smallest integer $k$ such that 
there exists a set $A$ of $k$ nonnegative integers such that $A$ is a basis of order $h$ for $n$, 
that is, $[0,n]\subseteq hA$.
There are, of course, many  open problems associated to the additive 
number theoretic  function $k_h(n)$. 

We have considered only sums of finite sets of integers and sums of finite sets of nonnegative 
integers, but there are analogous problems for sums of finite subsets of the squares, 
$k$th powers, primes, and other interesting sets of integers.  
Mathematicians have also considered not only the longest  intervals 
but also the longest finite arithmetic progressions, the longest finite geometric progressions,
and other finite configurations contained in $h$-fold sumsets.

Hans Rohrbach~\cite{rohr37a,rohr38} introduced the additive number theoretical functions 
$n_h(A)$ and $n_h(k)$ as follows. 
For every set $A$ of positive integers, $n_h(A)$ is the largest integer $n$ such that 
every integer in the interval $[1,n]$ is the sum of \emph{at most} $h$ not necessarily distinct 
elements of $A$.  
Equivalently, every every integer in the interval $[0,n]$ is the sum of \emph{exactly} $h$ 
not necessarily distinct elements of the set $A \cup \{0\}$ and so 
\[
n_h(A) = \ell_h(A \cup \{0\}). 
\]
For positive integers $h$ and $k$, Rohrbach defined $n_h(k)$ as the largest integer $n$ 
such that there exists a set $A$ of $k$ positive integers 
with $n_h(A) = n$ and so 
\[
n_h(k) = n^{\flat}_h(k+1).
\]

There is an extensive, albeit largely forgotten  literature, mostly in German, 
on the function $n^{\flat}_h(k)$ 
(e.g. Bloom~\cite{bloo26}, G{\"u}nt{\"u}rk and Nathanson~\cite{gunt-nath06}, 
H{\"a}mmerer and Hofmeister~\cite {hamm-hofm76},
H\" artter~\cite{hart60b}, Hofmeister~\cite{hofm63,hofm66a,hofm68,hofm83}, 
Kohonen~\cite{koho17}, 
 Marzuola and Miller~\cite{marz-mill10}, 
Moser~\cite{mose60}, 
Moser,  Pounder, and  Riddell~\cite{mose-poun-ridd69}, 
Mossige~\cite{moss81},
Mrose~\cite{mros79}, Nathanson~\cite{nath79b}, Rohrbach~\cite{rohr37a,rohr38}, 
Selmer~\cite{selm80}, 
St\" ohr~\cite{stoh55},  Weltge and Zyhalko~\cite{welt-zyha26}, Yu~\cite{yu09,yu15}) 
and it is possible 
that some of the  problems in this paper have solutions long  buried in back issues of journals, 
but now resident in and retrievable from LLMs.

\subsection*{Acknowedgements}
I thank Pietro Monticone and Wouter van Doorn 
for insightful discussions on the nature of AI problem solving in mathematics. 

\def\cprime{$'$} \def\cprime{$'$} \def\cprime{$'$}
\providecommand{\bysame}{\leavevmode\hbox to3em{\hrulefill}\thinspace}
\providecommand{\MR}{\relax\ifhmode\unskip\space\fi MR }
\providecommand{\MRhref}[2]{%
  \href{http://www.ams.org/mathscinet-getitem?mr=#1}{#2}
}
\providecommand{\href}[2]{#2}

\end{document}